\documentclass[12pt]{article}
\usepackage{amscd}
\usepackage{mathrsfs}
\usepackage{amsfonts}
\usepackage{amsmath,amscd,stmaryrd,amssymb,array, amsfonts,amsmath,amssymb,graphicx}
\usepackage{enumitem}
\setenumerate[1]{itemsep=0pt,partopsep=0pt,parsep=\parskip,topsep=5pt}
\setitemize[1]{itemsep=0pt,partopsep=0pt,parsep=\parskip,topsep=0pt}
\setdescription{itemsep=0pt,partopsep=0pt,parsep=\parskip,topsep=5pt}

\numberwithin{figure}{section}
 \numberwithin{equation}{section}

\newtheorem{theorem}{Theorem}[section]
\newtheorem{proposition}[theorem]{Proposition}
\newtheorem{definition}[theorem]{Definition}
\newtheorem{corollary}[theorem]{Corollary}
\newtheorem{lemma}[theorem]{Lemma}
\newtheorem{remark}[theorem]{Remark}

\newcommand{\cB}{{\mathcal B}}

\newcommand{\cD}{{\mathcal D}}

\newcommand{\cO}{{\mathcal O}}

\newcommand{\cK}{{\mathcal K}}

\newcommand{\cM}{{\mathcal M}}

\newcommand{\cV}{{\mathcal V}}

\newcommand{\cT}{{\mathcal T}}

\newcommand{\mB}{\mb{B}}

\newcommand{\sA}{{\mathscr A}}

\def\be{\begin{equation}}
\def\ee{\end{equation}}
\def\ba{\begin{array}}
\def\ea{\end{array}}
\def\bd{\begin{definition}}
\def\ed{\end{definition}}
\def\ben{\begin{enumerate}}
\def\een{\end{enumerate}}
\def\bt{\begin{theorem}}
\def\et{\end{theorem}}
\def\bp{\begin{proposition}}
\def\ep{\end{proposition}}
\def\bl{\begin{lemma}}
\def\el{\end{lemma}}
\def\br{\begin{remark}}
\def\er{\end{remark}}

\def\a{\alpha}

\def\de{\delta}
\def\pa{\partial}

\def\ve{\varepsilon}
\def\sig{\sigma}

\def\w{\omega}
\def\W{\Omega}
\def\gam{\gamma}

\def\.{\cdot}

\def\R{\mathbb{R}}

\def\A{\forall}
\def\ol{\overline}

\def\Cap{\bigcap}\def\Cup{\bigcup}

\def\ra{\rightarrow}

\def\~{\widetilde}
\def\8{\infty}
\def\X{\times}
\def\({\left(}
\def\){\right)}

\def\E{\exists}
\def\mb{\mbox}
\def\emp{\emptyset}
\def\-{\setminus}

\def\Hs{\hspace{0.8cm}}
\def\hs{\hspace{0.4cm}}
\def\Vs{\vskip8pt}
\def\vs{\vskip4pt}

\def\({\left(}\def\){\right)}

\begin{document}

\begin{center}
{\bf\Large  Attractors of Local Semiflows on \\[1ex] Topological Spaces 
}
\end{center}
\vskip10pt
\centerline{Desheng  Li\footnote[1]{Corresponding author. Supported by the grant of NSF of China (11071185, 11471240). {\em E-mail addresses}\,:  lidsmath@tju.edu.cn (D.S. Li).}
\hs Youbing  Xiong\,\, and  \,\, Jintao Wang
} 

\vskip20pt
\begin{center}
{\footnotesize
{Department of Mathematics, School of Science, Tianjin University\\
Center of Applied Mathematics, Tianjin University\\
     Tianjin 300072,  China}}
\end{center}

\vskip10pt

\begin{minipage}{13.5cm}
\centerline{\large\bf Abstract} \vskip10pt
In this  paper  we introduce a  notion of an attractor  for local semiflows on topological spaces, which in some cases seems to be more suitable  than the existing ones in the literature. Based on this notion   we   develop a basic attractor theory on topological spaces under appropriate separation axioms.  First, we discuss fundamental   properties of attractors such as maximality and stability and establish some existence results.  Then,  we give a converse Lyapunov theorem. Finally, the Morse decomposition of attractors is also addressed.
 \Vs

\Vs {\bf 2010 MSC:}  37B25, 35B34, 35B40, 35K55, 35J15.

\Vs\Vs {\bf Running Head:}  Attractors of Local Semiflows.


\end{minipage}

\newpage
\section{Introduction}
Invariant sets are of crucial importance  in the theory of dynamical systems.
This is because that for a given dynamical system, they are the  carriers of  much information on the longtime   behavior  of the system.   Of special interest are attractors. An
attractor, if exists,   is the depository of ``\,all\,'' the dynamics of a system near the attractor.

The attractor theories in metric spaces (especially
nonlocally compact metric spaces) were fully developed in the past
decades for both autonomous and nonautonomous systems \cite{Bab,Cheban,CV,Hale,Lady,Ma,Robin,Sell,Tem,Vishik}.
Here we are  interested in the case where the phase space is a  topological one that may not be metrizable. There are many motivations for this consideration. For example, for an infinite dimensional system on a Banach space $X$, in some cases one  has to study the dynamics of the system under weak topologies of $X$; see e.g. \cite{Gaz}. However, when  $X$ is endowed with its weak topology, it  may fail to be metrizable.  This makes the usual theory of dynamical systems in metric spaces inapplicable. Another example  is closely related to some recent  work  on the study of invariant sets. Let $N$ be an isolating neighborhood  of a semiflow $\Phi$ on a complete metric space $X$. Denote $N^-$ the {\em exit set} of $N$. Then one can define a quotient flow $\~\Phi$ on the quotient space $N/N^-$. It can be  shown under reasonable assumptions  that $\~\Phi$ is well-defined and continuous \cite{LiLK, WLD}. Thus many problems concerning invariant sets  (including the existence of invariant sets) of $\Phi$ can be transformed into that of attractors of the quotient flow $\~\Phi$.
But now a major obstacle  preventing us from taking a further step is that, in general  $N/N^-$ may not be metrizable, which makes the attractor theories on metric spaces fail to work.  To overcome this  difficulty,  a suitable attractor theory  on topological spaces needs to be developed.


In \cite{MN} Marzocchi and Necca introduced a notion of  {attractors } and established some existence results  for global semiflows on topological spaces. Informally speaking, the authors defined an attractor $\sA$ of a semiflow  on a topological space $X$ to be a compact invariant set  that attracts each element $B\in\cB$, where $\cB$ is a given family of subsets of $X$. Based on this notion  Giraldo et al. further  introduced the notion of  {global attractors} for global semiflows on topological spaces, using which they successfully  extended the shape theory of attractors from metric spaces to topological spaces \cite{GS}.

In practice, for an infinite dimensional  system $\Phi$, instead of imposing compactness conditions on the phase space $X$ one usually requires $\Phi$ to have  some {\em asymptotic compactness} properties. In the case where $X$ is non-metrizable, such a   property  seems to be more consistent with   {\em sequential compactness} of subsets of $X$. This simple observation stimulates us to introduce another notion of  attractors   for local semiflows on topological spaces. Specifically, we define an {attractor} to be a { sequentially compact invariant set} that  attracts a neighborhood of itself and enjoys some  maximal properties.  Since in a general topological space, compactness and sequential compactness are different matters, our notion of  attractors    differs significantly  from the ones in the literature, although we can show that they  coincide under appropriate conditions.

Based on our  notion of  attractors given here, we then  develop a basic attractor theory on topological spaces. First, we  discuss fundamental   properties of attractors such as maximality and stability. Then we give some  existence results. In particular, we show that if there is a closed set $M$  attracting an admissible neighborhood of itself, then $M$ contains an attractor.   Finally, the converse Lyapunov theorems and Morse decompositions of attractors are addressed. Our starting point is the convergence of sequences. It is worth mentioning  that throughout the paper we only assume  the phase space to be {\em  Hausdorff}  and {\em  normal}. No other separation axioms and countability axioms will be required. The interested reader is referred to \cite{LiLK} for an application of this attractor theory, in which we  proved some linking theorems and mountain pass type results to detect the existence of invariant sets of dynamical systems.

The theory of dynamical systems on topological spaces has a rich background. We refer  the  reader  to \cite{Bhatia, Got, Mark} etc. for some earlier work in this line.

This paper is organized as follows.
Section 2 is concerned with  fundamental properties of local semiflows, and Section 3 consists of some  results on limit sets. In Section 4 we  introduce the notion of an attractor, discuss basic properties of attractors and establish  some existence results. In Section 5 we prove a converse   Lyapunov theorem.  Section 6 is devoted to the Morse decomposition of attractors.

\section{Semiflows on topological spaces}


  Throughout the paper we  always assume that $X$ is a {\em Hausdorff}\, topological  space. Sometimes we may also require  $X$ to be $normal$, so that any two disjoint  closed subsets of $X$ can be separated  by their disjoint neighborhoods.

Let $A\subset X$. Denote $\ol A$, int$\,A$ and $\pa A$ the {\em closure, interior and {boundary}}  of any subset  $A$ of $X$,  respectively.
A set $U\subset X$ is called a {\em neighborhood} of $A$, if $\ol A\subset \mbox{int}\,U$.


We make a convention that we identify  a singleton $\{x\}$ with the point $x$.
\subsection{Definition and continuity property}
\begin{definition}\label{d2.1}\cite{Bhatia} A  local semiflow $\Phi$ on $X$ is a
continuous map from an open subset $\cD_\Phi$ of $\R^+\X X$  to $X$  satisfying the following conditions:
\begin{enumerate}
\item[$(1)$] For each $x\in X$, there exists $T_x\in(0,\8]$  such that $$(t,x)\in \cD_\Phi \Longleftrightarrow t\in[0,T_x)\,.$$
\item[$(2)$] $ \Phi(0,x)=x$ for all $x\in X$.
\item[$(3)$] If $(t+s,x)\in \cD_\Phi$, where $t,s\in\R^+$, then $$\Hs
\Phi(t+s,x)=\Phi\(t,\,\Phi(s,x)\).$$
\end{enumerate}

The set $\cD_\Phi$ and the number $T_x$ in the above definition are  called, respectively, the {domain} of $\Phi$  and the { escape time} of $\Phi(t,x)$.

A local semiflow $\Phi$ is called a {global semiflow}, if $\cD_\Phi=\R^+\X X$.
\end{definition}

Let $\Phi$ be a  given local semiflow on $X$. From now on  we rewrite $\Phi(t,x)$ as $\Phi(t)x$. For any sets $M\subset X$ and $J\subset \R^+$, denote
 $$
\Phi(J)M=\{\Phi(t)x:\,\,x\in M,\,\,t\in J\cap [0,T_x)\}.
$$

\bd
Given an interval $I\subset\R^1$, a map $\gamma:I\rightarrow X$ is called a {\em solution} (or {\em trajectory}) on
$I$, if $$\gamma(t)=
\Phi(t-s)\gamma(s),\Hs \A\,s,t\in I,\,\,s\leq t.$$
A  solution $\gam$ on $I=\R^1$ is  called a {\em full solution}.
\ed

 It is known that a solution is continuous \cite{Bhatia}.
Let $\gam$ be a solution  on $I$. Set $$\cT_\gamma(I)=\{(t,\gamma(t)):\,\,t\in I\}.$$   $\cT_\gam(I)$ is called the {\em trace} of $\gam$ on $I$.

For any $x\in X$ and $I\subset \R^+$, we will write $\cT_x(I)=\cT_{\gam_x}(I)$, where $\gam_x(t)=\Phi(t)x$.


\bp\label{l2.3} Let $x\in X$, and  $0\leq T<T_x$. Then \begin{enumerate}\item[$(1)$] there exists a neighborhood $U_x$ of $x$ such that $$T_y>T,\Hs \A\,y\in U_x;$$
\vs \item[$(2)$] for any compact interval $J\subset [0,T]$, $\cT_{y}(J)\ra \cT_{x}(J)$ in $J\X X$ as $y\ra x$. Specifically,
for any neighborhood $\cV$ of \,$\cT_{x}(J)$ in $J\X X,$   there exists a neighborhood $U$ of $x$ such that
\be\label{e2.0}
\cT_{y}(J)\subset \cV,\Hs \A\,y\in U.\ee
\end{enumerate}
\ep

\noindent{\bf Proof.}  The first conclusion (1) is implied in \cite{Bhatia}, Lemma 1.8.
To prove the second one, we define a map $G:\cD_\Phi\ra\R^+\X X$ as follows:
$$
G(t,x)=(t,\Phi(t)x),\Hs (t,x)\in\cD_\Phi.
$$
Clearly $G$ is continuous. Let $U_x$ be the neighborhood of $x$ in (1), and $\cV$ be a neighborhood  of \,$\cT_{x}(J)$ in $J\X X$. Then for each $t\in J$, there exists a cylindrical neighborhood $Q_t:=I_t\X U_t$ of $(t,x)$ such that $G(Q_t)\subset \cV$, where $I_t$ is an interval relatively open in $J$, and $U_t\subset U_x$ is a neighborhood of $x$. Since $J\subset \Cup_{t\in J}I_t$,  there exists a finite number of ${I_t}'s$, say, $I_{t_1}, I_{t_2},\cdots, I_{t_n}$, such that $J\subset \Cup_{i=1}^nI_{t_i}$. Set $U=\Cap_{i=1}^n U_{t_i}$. Then $U$ is a neighborhood of $x$. It can be seen that $U$ fulfills \eqref{e2.0}.
\,$\Box$
\br If $X$ is a  metric space with metric $d$, then  Pro.\,\ref{l2.3} can be reformulated in a simpler but more specific  manner as below:
\bp\label{p:2.5} Let $x\in X$, and $0<T<T_x$. Then for any $\ve>0$, there exists  $\de>0$ such that $\Phi(t)y$ exists on $[0,T]$ for all $y\in \mB(x,\de)$. Moreover,
$$
d\(\Phi(t)y,\,\Phi(t)x\)<\ve,\Hs \,\A\,t\in[0,T],\,\,y\in \mB(x,\de).
$$
 \ep
 \er


\subsection{A convergence result of solutions}
In this subsection we give a convergence result concerning sequences of solutions.

Let us first  recall the concepts of  {\em strong admissibility}, which was first introduced for local semiflows on metric spaces \cite{Ryba}.
%

Let $M$ be a subset of $X$.

\begin{definition}\label{d4.6c}
We say that  $\Phi$ does not explode in $M$, if  $T_x=\8$ whenever $\Phi([0,T_x))x\subset M.$
\end{definition}

\begin{definition}\label{d4.6} $M$  is said to be  {admissible}, if for any sequences $x_n\in M$ and $t_n\ra \8$ with $\Phi([0,t_n])x_n\subset M$ for all $n$,  the sequence $\Phi(t_n)x_n$ has a convergent subsequence.

$M$ is said to be strongly admissible, if it is admissible and moreover, $\Phi$ does not explode in $M$.
\end{definition}

\bp\label{l2.9} Let $M$ be a closed  strongly admissible set, and  $\gamma_n$ be a sequence of solutions in $M$ with each  $\gamma_n$ being defined on $[-t_n,t_n]$.

Suppose  $t_n\ra\8$. Then  there exists a  subsequence $\gamma_{n_k}$ of $\gamma_n$ and a full solution $\gamma$  in $M$, such that for any compact interval $J\subset \R^1$, $$\mb{$\cT_{\gamma_{n_k}}(J)\ra\cT_\gamma(J)$ in $J\X X$.}$$ That is,
for any neighborhood $\cV$ of \,$\cT_\gamma(J)$ in $J\X X,$   there exists $k_0>0$ such that
\be\label{e2.6}
\cT_{\gamma_{n_k}}(J) \subset \cV,\Hs \A\,k>k_0.\ee
\ep
\noindent{\bf Proof.} We may assume $t_n>1$ for all $n$. Since $\gamma_n(-1)=\Phi(t_n-1)\gamma(-t_n)$, by  admissibility of $M$    there exists a subsequence $\gamma_{1n}$ of $\gamma_n$ such that $\gamma_{1n}(-1)$ converges to a point $x_1\in M$. Set
$$\sig_1(t)=\Phi(t+1)x_1,\Hs t\in[-1,T_{x_1}-1).$$
 We claim that $\sig_1$ is a solution  on $[-1,\8)$ contained in $M$.

 Indeed, since $\gamma_{1n}(-1)\ra x_1$, by continuity of $\Phi$ we deduce that $\gamma_{1n}(t)\ra \sig_1(t)$ for all $t\in [-1,T_{x_1}-1)$.
The  closedness of $M$ then implies that $\sig_1(t)\in M$ for all $t\in [-1,T_{x_1}-1)$, i.e.,
 $\Phi(t)x_1\in M$ for all $t\in[0,T_{x_1})$. As $\Phi$ does not explode in $M$, we have $T_{x_1}=\8$. Thus the claim holds true.

We also infer from  Pro. \ref{l2.3}  that $\cT_{\gamma_{1n}}(J)$ converges to $\cT_{\sig_1}(J)$ in $J\X X$ for any compact interval $J\subset [-1,\8)$.

 Repeating the same argument with very minor modifications, one can show that there exist a subsequence $\gamma_{2n}$ of $\gamma_{1n}$ and a solution $\sig_2$  on $[-2,\8)$ contained in $M$ such that  $\cT_{\gamma_{2n}}(J)$ converges to $\cT_{\sig_2}(J)$ in $J\X X$ for any  compact interval $J\subset [-2,\8)$. Clearly $$\sig_2(t)\equiv\sig_1(t),\Hs t\in [-1,\8).$$

 Continuing the above  procedure, we obtain for each $k$ a subsequence $\gamma_{kn}$ of $\gamma_n$ and a solution $\sig_k$  in $M$ such  that
 \begin{enumerate}
 \item[(1)] each sequence $\gamma_{(k+1)n}$ is a subsequence of $\gamma_{kn}$;\vskip-5pt
 \item[(2)] $\sig_k$ is defined on $[-k,\8)$, and \be\label{e2.3}\sig_{k+1}(t)\equiv\sig_k(t),\Hs t\in [-k,\8);\ee
 \item[(3)] for any compact interval $J\subset [-k,\8)$, $\cT_{\gamma_{kn}}(J)$ converges to $\cT_{\sig_k}(J)$ in $J\X X$ as $n\ra\8$ in the sense that for any neighborhood $\cV$ of \,$\cT_{\sig_k}(J)$ in $J\X X,$   there exists $n_k>0$ such that
$$
\cT_{\gamma_{kn}}(J) \subset \cV,\Hs \A\,n>n_k.$$
 \end{enumerate}
 Define a full solution $\gamma$ in $M$ as follows:
 $$
 \gamma(t)=\sig_k(t),\Hs \mb{if }\,t\in[-k,\8).
 $$
By \eqref{e2.3} it is clear that  $\gamma$ is well defined. Consider the sequence $\gamma_{kk}$\,. By virtue of the classical  diagonal procedure it can be easily seen that  $\cT_{\gamma_{kk}}(J)$ converges to $\cT_{\gamma}(J)$ in $J\X X$ as $k\ra \8$ for any compact interval $J\subset \R^1$. \,$\Box$

\section{Invariant Sets and Limit sets} In this section we talk about some basic facts on invariant sets and  limit sets of local semiflows on topological spaces.

Let $X$ be a Hausdorff topological space, and $\Phi$ be a local semiflow on $X$.

\bd
The {\em $\omega$-limit set} $\w(M)$ of $M\subset X$ is defined as
$$\ba{ll}
\w(M)=\{y\in X:\,\,\,\,\E\, x_n\in M\mbox{ and }
t_n\rightarrow+\infty \mb{ such that } \Phi(t_n)x_n\rightarrow y\}.
\ea
$$

The {\em $\omega$-limit set} $\w(\gamma)$ of a solution $\gam$ on $(a,\8)$ is defined as
$$\ba{ll}
\w(\gamma)=\{y\in X:\,\,\,\,\E\,  t_n\ra \8 \mb{ such that }\gamma(t_n)\ra y\}.\ea$$
The {\em $\alpha$-limit set} $\alpha(\gamma)$ of a solution $\gam$ on $(-\8,a)$ is defined as
$$\ba{ll}
\alpha(\gamma)=\{y\in X:\,\,\,\,\E\,  t_n\ra -\8 \mb{ such that }\gamma(t_n)\ra y\}.\ea$$
\ed
\br In \cite{MN} the $\w$-limit set $\w(M)$ of $M\subset X$ is defined as
$$\w(M):=\Cap_{t\geq 0}\ol{\Phi([t,\8))M},$$
which is somewhat different from that of ours here. But if we comeback to the situation of a metric space, then both  coincide \cite{Tem}.
\er

To discuss fundamental  properties of limit sets, we need to introduce several  notions on invariance and attraction.

\vs Let  $M$ be a subset of $X$.
\bd
$M$ is said to be  { positively (resp.  negatively) invariant}\, if
$$\Phi(t)M\subset M\,\,\,(\mb{resp. }M\subset \Phi(t)M),\Hs \A\,t\geq 0.$$
$M$ is said to be { invariant}, if it is both negatively and positively invariant.

\ed

\bp\label{p3.4}\cite{Bhatia} If $M$ is positively invariant, then so is $\ol M$.
\ep



\br We do not know whether a similar result holds true for negative invariance.
  But by a very standard argument one can easily verify that if $M$ is negatively invariant, then  for each $y\in M$ there exists a solution $\gamma$ on $(-\8,0]$ contained  in $M$ such that $\gamma(0)=y$.
Further noting that $y=\Phi(n)\gam(-n)$ for all $n$, by the definition of $\w(M)$ one concludes  that $y\in \w(M)$.
  Hence
    $$M\subset \w(M).$$
 \er


\begin{definition} We say that  $M$ {\em attracts} $B\subset X$, if
 $T_x=\8$ for all $x\in B$ and moreover,
 for any neighborhood $V$ of $M$, there exists $t_0>0$ such that $$\Phi(t)B\subset V,\Hs\A\,t>t_0.$$
\end{definition}

\begin{proposition}\label{p2.2} Suppose  ${\Phi(\R^+)M}$  is contained in a closed  strongly admissible set $N$.
 Then $\w(M)$ is a nonempty invariant set  that attracts $M$.
\end{proposition}

\noindent{\bf Proof.} The  nonemptyness of  $\w(M)$ is a simple consequence of the strong admissibility of $N$. Clearly $\w(M)\subset N$.
We show that $\w(M)$ is  invariant.

 Let $y\in \w(M)$. Then there exist sequences $x_n\in M$ and $t_n\ra\8$ such that $y_n:=\Phi(t_n)x_n\ra y$. We claim that  $\Phi(t)y\in\w(M)$ for all $t\in[0,T_y)$, hence $\w(M)$ is positively invariant.
 Indeed, for any $0<t<T_y$ we have
$$
\Phi(t)y=\lim_{n\ra\8}\Phi(t)y_n= \lim_{n\ra\8}\Phi(t)\Phi(t_n)x_n=\lim_{n\ra\8}\Phi(t+t_n)x_n,
$$
from which one immediately deduces that  $\Phi(t)y\in \w(M)$. 


We also show that for any $t\geq 0$,  there exists $z\in \w(M)$ such that $\Phi(t)z=y$, thus $\w(M)$ is negatively invariant.
We may assume $t_n\geq t$ for all $n$.
 Let $z_n=\Phi(t_n-t)x_n$. By  admissibility of $N$ there exists a subsequence of $z_n$ (still denoted by $z_n$) converging  to a point $z$. Clearly  $z\in \w(M)$.  We have
$$
y=\lim_{n\ra\8}\Phi(t_n)x_n=\lim_{n\ra\8}\Phi(t)\Phi(t_n-t)x_n=\lim_{n\ra\8}\Phi(t)z_n=\Phi(t)z.
$$

The verification of the attraction property   is trivial. We omit the details.
\, $\Box$
\Vs

\begin{proposition}\label{p2.1} Let $\gamma$ be a solution on $I:=(a,\8)$ (resp. $(-\8,a))$). Suppose  $\gamma(I)$ is contained in a closed strongly  admissible set $N$. Then
$\w(\gam)$ (resp. $\a(\gamma)$\,) is a nonempty invariant set.
\end{proposition}

\noindent{\bf Proof.} If $\gamma$ is  a solution on $(a,\8)$, then it can be easily seen that  $\w(\gam)=\w(\gamma(t_0))$, where  $t_0>a$. Thus the conclusion on $\w(\gam)$  immediately follows from Pro.\,\ref{p2.2}.

Now let  $\gamma$ be  a solution on $(a,\8)$ and consider the $\alpha$-limit set $\alpha(\gam)$.  For  convenience in statement, we may assume $a=0$.
 We first check that $\alpha(\gam)\ne\emp$. Take a sequence $0<t_n\ra \8$. Then $\gam(-t_n)=\Phi(t_n)\gam(-2t_n)$. By  admissibility of $N$ we deduce  that  $\gam(-t_n)$ has a  subsequence converging to a point $z$. Clearly  $z\in \alpha(\gam)$.

As $N$ is closed, we have $\alpha(\gamma)\subset N$.
 Let $y\in\alpha(\gam)$, and let $t_n\ra\8$ be such that $y_n:=\gam(-t_n)\ra y$.
By  a  similar argument as in the  verification of  positive invariance of $\w(M)$ in Pro.\,\ref{p2.2}, it can be shown that $\Phi(t)y\in\alpha(\gamma)$ for all $t\in[0,T_y)$. 
Therefore  $\alpha(\gamma)$ is  positively invariant.
The verification of negative invariance of $\alpha(\gam)$ can also be performed in a similar manner as in the case of $\w(M)$ in Pro.\,\ref{p2.2}, and is thus omitted.
\,$\Box$
\section{Attractors}
In this section we introduce the notion of attractors and discuss basic properties of attractors. Some existence results will also be given.
\vs
Let $X$ be a Hausdorff space, and $\Phi$ be a given local semiflow on $X$.
\subsection{Definition and basic properties}
We first recall that a set $M\subset X$ is said to be {\em sequentially compact} ({\em s-compact} in short), if each sequence $x_n$ in $M$ has a subsequence converging to a point $x\in M$.

\begin{definition} \,A nonempty s-compact invariant set $\sA$ is called an attractor of $\Phi$,  if
there is a neighborhood $N$ of $\sA$ such that
\ben
\item[$(1)$] $\sA$ attracts $N$; and
\item[$(2)$]  $\sA$ is the maximal s-compact invariant set in $N$.
\een
\ed

Given  an attractor $\sA$, define $$ \W(\sA)=\{x\in X:\,\,\sA\mb{ attracts }x\}.
$$
 $\W(\sA)$ is called the {\em region of attraction } of $\sA$. If $\W(\sA)=X$, then  $\sA$ is simply called a {\em global attractor}


\bt\label{p2.3} Let $\sA$ be an attractor of $\Phi$. Then the following assertions hold.
\begin{enumerate}
\item[$(1)$] $\W(\sA)$ is open, and for each compact set $K\subset \W(\sA)$,  $\sA$ attracts a   neighborhood $U$ of $K$.\vs
\item[$(2)$] $\sA$ is the maximal s-compact invariant set  in $\W(\sA)$.\vs
\item[$(3)$] If $X$ is normal, then for any closed admissible neighborhood $V$ of $\sA$ with $V\subset \W(\sA)$, $\sA$ is the maximal invariant set in $V$.
\end{enumerate}
\et

\noindent{\bf Proof.} (1) By  definition  there is a neighborhood $N$ of $\sA$ such that $\sA$  attracts $N$, moreover,   $\sA$ is the maximal s-compact invariant set in $N$.

Take a $\tau>0$ such that
\be\label{e2.4}
\Phi(t)N\subset \mb{int}N,\Hs t\geq \tau.
\ee

Let $x\in\W(\sA)$. Then the escape time  $T_x=\8$. Furthermore,  there exists $T>0$ such that  $$\Phi(T)x\in\mb{int}N.$$ By virtue of Pro.\,\ref{l2.3}, there is a neighborhood $U_x$ of $x$ such that $T_y>T$ and $\Phi(T)y\in \mb{int}N$ for all $y\in U_x$, from which it can be easily seen that  $\sA$ attracts each point in $U_x$. Hence  $\W(\sA)$ is  open.

\vs Let $K$ be a compact subset of $\W(\sA)$. 
 For each $x\in K$, pick a $t_x>0$ such that $\Phi(t_x)x\in \mb{int}N$. Then by continuity  there exists an open  neighborhood $U_x$ of $x$ such that $\Phi(t_x)U_x\subset \mb{int}N$. Combining this with \eqref{e2.4} it yields
\be\label{e2.5}
\Phi(t)U_x\subset \mb{int}N,\Hs t\geq \tau+t_x:=\tau_x.
\ee
Since  $K$ is compact,  there exists a finite number of points in $K$, say, $x_1,x_2,\cdots,x_n$ such that $$K\subset \Cup_{1\leq i\leq n}U_{x_i}:=U.$$ Let $\tau_0=\max\{\tau_{x_i}:\,\,{1\leq i\leq n}\}$.
Then \eqref{e2.5} implies  that $\Phi(\tau_0)U\subset N$. Hence we see that   $\sA$ attracts $U$.

\vs(2) Let $M$ be an s-compact invariant set in $\W(\sA)$. We need to prove that $M\subset \sA$.
Since $\sA$ is the maximal s-compact invariant set of $\Phi$ in $N$, for this purpose it suffices to check that $M\subset N$.

We argue by contradiction and suppose the contrary.  Then there would exist $y\in M\-N$.
Let $\gam$ be a solution on $(-\8,0]$ contained in $M$ with $\gam(0)=y$. As $M$ is   s-compact, using some similar argument as in the proof of Pro.\,\ref{p2.1} one can easily verify  that $\alpha(\gam)$ is a nonempty invariant set with $\alpha(\gam)\subset M\subset \W(\sA)$.

Now if $\alpha(\gam)\cap \,\mb{int}N\ne\emp$,  one can pick a  $z\in\alpha(\gam)\cap \,\mb{int}N$.  Take   a sequence $t_n\ra\8$ such that $\gam(-t_n)\ra z$. Then $\gam(-t_n)\in \mb{int}N$ for $t_n$ sufficiently large. Fix a  $t_n>\tau$ such that $\gam(-t_n)\in \mb{int}N$. By \eqref{e2.4} we find that
$$y=\gam(0)=\Phi(t_n)\gam(-t_n)\in N,$$ a contradiction!
On the other hand, if  $\ba{ll}\alpha(\gam)\cap \mb{int}N=\emp,\ea$
then by invariance of $\alpha(\gam)$ we necessarily have  $\alpha(\gam)\cap\, \W(\sA)=\emp$, which again leads to a contradiction.

\vs
(3) Assume $X$ is normal. Let  $V\subset \W(\sA)$ be a closed  admissible neighborhood of $\sA$. Then one can find a closed neighborhood $O$ of $\sA$ with $O\subset V\cap N$. Clearly $\sA$ attracts $O$. Hence there exists $t_0>0$ such that \be\label{e4.9}\Phi([t_0,\8))O\subset O.\ee By \eqref{e4.9} and  Pro.\,\ref{p2.2} we deduce that $\w(O)$ is a nonempty invariant set with $\w(O)\subset O\subset N$. We claim that
$\w(O)$ is s-compact.
Indeed, let $y_n\in \w(O)$ be a sequence. Then by invariance of $\w(O)$ there exists a sequence $x_n\in\w(O)$ such that $y_n=\Phi(n)x_n$. The admissibility of $O$ implies that $y_n$ has a subsequence $y_{n_k}$ converging to a point $y$. As $x_n\in\w(O)\subset O$, by the definition of $w$-limit set we have $y\in\w(O)$, which completes the proof of the claim.

Because $\sA$ is the maximal $s$-compact invariant set in $N$ and $\w(O)\subset N$, we have \be\label{e4.10}\w(O)\subset \sA.\ee
Now if  $K$ is  an invariant set in $O$, then $K\subset \w(O)$. Hence  by \eqref{e4.10} we have $K\subset \sA$. It follows that  $\sA$ is the maximal invariant set in $O$.

To show  that $\sA$ is the maximal invariant set in $V$, it now suffices  to check that if $K$ is an invariant set in $V$, then $K\subset O$.
We argue by contradiction and suppose $K\-O\ne\emp$. Pick a $y\in K\-O$. Then there is a solution $\gam$ on $(-\8,0]$ contained in $K$ such that  $\gam(0)=y$.
Since  $K\subset V$, by Pro.\,\ref{p2.1} $\alpha(\gam)$ is a nonempty invariant set.  Using some similar argument as in the proof of (2) (with $N$ therein replaced by $O$), one immediately obtains a contradiction. \, $\Box$

\br We infer from  Theorem \ref{p2.3} (1) that  the global attractor, if exists, is necessarily unique.

In \cite{GS}  a global attractor is defined  to be a compact invariant set that attracts each compact set. By Theorem\,\ref{p2.3} (1) we see that if a global attractor  in our terminology is compact, then it is
a global attractor in the terminology in \cite{GS}.
\er

Now we turn our attention to stability of attractors.

\begin{definition}
 A set $M\subset X$ is called  {stable}, if
for any neighborhood $U$ of $M$, there is a neighborhood $V$ of $M$
such that \be\label{e4.4}\Phi(\R^+)V\subset U.\ee
\end{definition}

\bt\label{st} Let $\sA$ be an attractor of $\Phi$. Then $\sA$ is stable.\et

\noindent{\bf Proof.} Let $U$ be an open neighborhood of $\sA$.  We  need to prove that there exists an open  neighborhood $V$ of $\sA$ such that
\ref{e4.4} holds.

By  definition   $\sA$ attracts a neighborhood $N$ of $\sA$.
Fix a $T>0$ so that
\be\label{e4.6}\Phi(t)N\subset U,\Hs \A\,t>T.\ee
 In what follows  we argue  by contradiction and suppose \eqref{e4.4} fails to be true. Then for any open neighborhood $V\subset N$ of $\sA$, there exist $y\in V$ and $t_y\geq0$ such that
\be\label{e4.5}\Phi(t_y)y\notin U.\ee In view of \eqref{e4.6}, we necessarily  have $t_y\leq T$.

\vs We claim that  there exists $x_0\in\ol\sA$ such that for any open neighborhood $O\subset N$ of $x_0$, one can find a $z\in O$ and $t_z\in[0,T]$ such that
\be\label{e4.2}
\Phi(t_z)z\notin U.\ee
Indeed, if this was false,  each point $x\in\ol\sA$ would have  an open neighborhood $U_x\subset N$ such that $\Phi([0,T])U_x\subset U$.  Set $V=\Cup_{x\in\ol\sA}U_x$. Obviously  $V$ is an open  neighborhood of $\sA$. However, $\Phi([0,T])V\subset U$, which  contradicts \eqref{e4.5} and proves our claim.

On the other hand,  we infer from   Pro.\,\ref{p3.4} that $\ol\sA$ is positively invariant.
Hence $\Phi([0,\,T])x_0\subset\ol\sA\subset U$. As $U$ is open, by virtue of Pro.\,\ref{l2.3}   there is an open neighborhood $O$ of $x_0$ such that $\Phi([0,\,T])O\subset U$. But this contradicts \eqref{e4.2}.

 The proof of the theorem is finished. \,$\Box$

\subsection{Existence  of attractors}
 Now we  prove some existence results on attractors. First, we have

\bt\label{p2.4}  Suppose $X$ is normal. Let $M\subset X$ be closed. Assume $M$ attracts an  admissible
neighborhood $N$ of itself.  Then  $\sA=\w(N)$ is an attractor.
\et
\noindent{\bf Proof.} Since $X$ is normal, it can be assumed that $N$ is closed.  By the  definition of attraction we have $T_x=\8$ for all $x\in N$. Furthermore,  there exists $T>0$ such that $$\Phi([T,\8))N\subset N.$$ Pro.\,\ref{p2.2} then asserts   that $\sA=\w(N)$ is a nonempty invariant set  attracting $N$.

We claim that $\sA\subset M$. Indeed, if $y\in \sA\- M\ne \emp$, then by  closedness of  $M$ there exist open neighborhoods $U$ of $M$ and $V$ of $y$ such that $U\cap V=\emp$. Let $x_n\in N$ and $t_n\ra\8$ be such that $\Phi(t_n)x_n\ra y$. Then $\Phi(t_n)x_n\in V$ for $n$ sufficiently large. On the other hand, since $M$ attracts $N$, one should have
$\Phi(t_n)x_n\in U$ for $n$ sufficiently large, which leads to a contradiction. Hence the claim holds true.

  Let $K\subset N$ be an invariant set. Then
$$
K\subset \w(K)\subset \w(N)=\sA.
$$
Thus we deduce that $\sA$ is the maximal invariant set in $N$.

In what follows we show that  $\sA$ is s-compact,  hence it is an attractor.
Let $y_n$ be a sequence in $\sA$. Then  for each $n$ there is an $x_n\in \sA$ such that $y_n=\Phi(n)x_n$. Further by  admissibility of $N$ we deduce that $y_n$ has a convergent subsequence $y_{n_k}$. Let $y_{n_k}\ra y$ as $k\ra\8$. We observe that  $$y=\lim_{k\ra\8}y_{n_k}=\lim_{k\ra\8}\Phi(n_k)x_{n_k}\in\w(N)= \sA.$$
Therefore one concludes that   $\sA$ is s-compact. \,$\Box$


\bt\label{p2.5} Suppose $X$ is normal. Let $M$ be a closed subset of $X$. Assume $M$ is stable, and  that there exists a strongly  admissible  neighborhood $N$ of $M$  such that  $M$ attracts each point  $x\in N$.

 Then $\Phi$ has  an attractor $\sA\subset M$ with $N\subset \W(\sA)$.
\et

\noindent{\bf Proof.} Take a closed neighborhood $F$ of $M$ with $F\subset N$. Then by stability of  $M$   there is  a closed neighborhood $W$ of $M$ such that $\Phi([0,\8))W\subset F$.
 Thus by Pro.\,\ref{p2.2} we deduce  that   $\w(W)$ is an invariant set that attracts $W$. Clearly $\w(W)\subset F$.
We  show that
\be\label{e4.7}\w(W)\subset M,\ee
hence $M$ attracts $W$. It then  immediately follows by Theorem \ref{p2.4} that $\sA=\w(W)$ is an attractor.

 We argue by contradiction and suppose the contrary.  Then $\w(W)\-M\ne\emp.$ Pick a $y\in \w(W)\-M$. One can find  a neighborhood $V$ of $M$ such that $y\not\in V$.
Take  an open  neighborhood $U$ of $M$ such that $$\Phi(t)U\subset V,\Hs \A\,t\geq 0.$$
Let $\gam$ be a solution on $(-\8,0]$ contained in $\w(W)$ with $\gam(0)=y$. If $\gam(\tau)\in U$ for some $\tau<0$, then $$y=\gam(0)=\Phi(-\tau)\gam(\tau)\in V,$$
which leads to a contradiction. Thus we have \be\label{e4.11}\gam(t)\not\in U,\Hs \A\,t<0.\ee
Noticing that  $\gam((-\8,0])\subset \w(W)\subset F$, by Pro.\,\ref{p2.1} we see that $\alpha(\gam)$ is a nonempty  invariant set. Clearly $\alpha(\gam)\subset F$.
We claim that \be\label{e4.8}\alpha(\gam)\cap U=\emp.\ee Indeed, if this was false, there would  exist $z\in \alpha(\gam)\cap U$. Take a   sequence $t_n\ra -\8$ such that $\gam(t_n)\ra z$. Then $\gam(t_n)\in U$ for $t_n$ sufficiently large,
which  contradicts \eqref{e4.11} and proves \eqref{e4.8}.

Now since  $\alpha(\gam)$ is invariant,   \eqref{e4.8} implies that $M$ does not attract any point $x\in\alpha(\gam)$, which contradicts the attraction assumption on $M$ (recall that $\alpha(\gam)\subset F\subset N$) and completes the proof of \eqref{e4.7}.

To complete the proof of the theorem, there remains to check that  $N\subset \W(\sA)$. Let $x\in N$. As  $M$ attracts $x$ and $W$ is a neighborhood of $M$,  there exists $t_0$ such that $\Phi(t_0)x\in W$. As  $\sA$ attracts $W$, it immediately follows that  $x\in\W(\sA)$. \, $\Box$
\Vs

The following result is a simple  consequence of Theorem \ref{p2.5}. It can be seen as a converse theorem of Theorem \ref{st}.
\begin{corollary} Let $\sA$ be  a closed invariant set. If  $\sA$ is stable and attracts each point in a strongly admissible  neighborhood of itself, then it is an attractor.
\end{corollary}

\section{Lyapunov Functions of Attractors}
It is well known that Lyapunov functions play crucial roles in many aspects of stability analysis.
In this section we prove a converse  Lyapunov theorem for attractors of local semiflows on topological spaces.

Let  $\Phi$ be a given local semiflow on Hausdorff space $X$, and  $\sA$ be an attractor of $\Phi$ with the region of attraction  $\W=\W(\sA)$.
\bd
A nonnegative  function $\zeta\in C(\W)$ is called a $\cK_0$ {\em function} of $\sA$, if
$$
\zeta(x)=0\Longleftrightarrow x\in\sA.
$$
\ed
\begin{definition}
 A $\cK_0$ {function} $L$ of $\sA$ is called a Lyapunov function of $\sA$, if
$$
L(\Phi(t)x)<L(x),\Hs\A\ x\in \W\-\sA,\,\,t>0.
$$
\end{definition}



\bt\label{p2.13}Assume that  $\sA$ is closed and  has a $\cK_0$ function $\zeta$  on $\W$.
Then $\sA$ has a Lyapunov $L$ on $\W$.

If we assume, in addition, that  $X$ is normal, then  for any  closed subset $K$ of $X$ with  $K\subset \W\- \sA$, there exists  a Lyapunov function $L$ of $\sA$ such that
\be\label{cl}
L(x)\geq 1,\Hs \A\,x\in K.
\ee
\et
\noindent{\bf Proof.}
The  construction of $L$ is the same as  in \cite{Kap} (pp. 226). See also \cite{Li3}. We give the details for the reader's convenience.

 Set
\be\label{n2}\xi(x)=\sup_{t\geq0} \zeta(\Phi(t)x),\Hs x\in \W.\ee
We show that $\xi$ is continuous. For this purpose, it suffices to check that
for any fixed $x_0\in\W$ and $\ve>0$,
there is a neighborhood $U$ of $x_0$ with $U\subset \W$ such that
\be\label{e5.3}|\xi(x)-\xi(x_0)|<\ve,\Hs\A\, x\in U.\ee

Let $\xi_0=\xi(x_0)$. First, by the definition of $\xi$ there exists $t_0\geq 0$ such that $$\zeta(\Phi(t_0)x_0)>\xi_0-\ve/2.$$ Hence by  continuity  one can find
a neighborhood $V$ of $x$ with $V\subset \W$ such that
$$\zeta(\Phi(t_0)x)>\xi_0-\ve,\Hs\A\,x\in V.$$ Therefore
\be\label{e5.4}
\xi(x)\geq \zeta(\Phi(t_0)x)>\xi_0-\ve,\Hs\A\,x\in V.\ee
We also check that there is a neighborhood $U$ of $x_0$ with $U\subset V$ such that
\be\label{e5.5}\xi(x)<\xi_0+\ve,\Hs \A\,x\in U,\ee
thus proving  \eqref{e5.3}.

Let $W=\zeta^{-1}([0,\,\xi_0+\ve/2))$. Then $W$ is an open neighborhood of $\sA$. By  definition   $\sA$ attracts a  neighborhood $N$ of itself. It can be assumed that $N$ is open.
Fix a  $T_1>0$ such that \be\label{e5.6}\Phi(t)N\subset W,\Hs\A\,t\geq T_1.\ee
Take a $T_2>0$ such that
 $\Phi(T_2)x_0\in N$. Then by continuity of $\Phi$ one can find  a neighborhood $U'$ of $x_0$ with  $U'\subset V$ such that $\Phi(T_2)U'\subset N$.
Combining this with \eqref{e5.6} it yields
$$\Phi(t)U'\subset W,\Hs \A\,t>T_1+T_2:=T.$$
Hence by the definition of $W$ we have
\be\label{n3}\zeta(\Phi(t)x)< \xi_0+\ve/2,\Hs\A\,x\in U',\,\,t>T.\ee

Recalling that $\zeta(\Phi(t)x_0)\leq\xi(x_0)=\xi_0$ for $t\geq 0$, we have $$\Phi([0,T])x_0\subset W.$$  Thus by Pro.\,\ref{l2.3} one easily deduces that there is a neighborhood $U$ of $x_0$ with  $U\subset U'\subset V$ such that
$\Phi([0,\,T])U\subset W.$ It follows that
\be\label{n4}\zeta(\Phi(t)x)<\xi_0+\ve/2,\Hs \A\,x\in U,\,\,t\in[0,\,T].\ee
(\ref{n3}) and (\ref{n4}) assures that $$\zeta(\Phi(t)x)<\xi_0+\ve/2,\Hs\,\A x\in U,\,\,t\geq 0.$$ Thereby
 $$\xi(x)=\sup_{t\geq0} \zeta(\Phi(t)x)\leq\xi_0+\ve/2<\xi_0+\ve,\Hs \A\,x\in U.$$ This  completes the proof of \eqref{e5.5}.
\Vs
Clearly $\xi(x)\equiv 0$ on $\sA$. A basic property of $\xi$ is that it is decreasing along each solution of $\Phi$ in $\W$.
Indeed, for any $x\in\W$ and $t\geq 0$, we have $$ \xi(\Phi(t)x)=\sup_{s\geq0} \zeta\(\Phi(s)\Phi(t)x\)=\sup_{\tau\geq t}\zeta(\Phi(\tau)x)\leq \xi(x)$$

Define $$ L(x)=\xi(x)+\int_0^\8e^{-t}\xi(\Phi(t)x)dt,\Hs x\in \W. $$ We show that $L$ is precisely a Lyapunov function of $\sA$.

Let $x\in\W\setminus\sA$, and  $s>0$. Then $\xi(x)>0$.
Because  $\xi$ is decreasing along each solution  in $\W$, we see  that
$$
\xi(\Phi(t)\Phi(s)x)=\xi\(\Phi(s)\Phi(t)x\)\leq \xi\(\Phi(t)x\),\Hs \A t\geq 0.
$$
We claim that there is at least one point $t\geq 0$ such that $$\xi\(\Phi(t)\Phi(s)x\)< \xi\(\Phi(t)x\).$$ Indeed, if this was false,  one should have
$$
\xi\(\Phi(t)x\)\equiv\xi\(\Phi(t)\Phi(s)x\)=\xi\(\Phi(t+s)x\),\Hs t\geq 0.
$$
Hence $\xi\(\Phi(t)x\)$ is an $s$-periodic function. This contradicts to the fact that $\xi\(\Phi(t)x\)\ra0$ as $t\ra\8$.

Now since  both $\xi\(\Phi(t)\Phi(s)x\)$ and $\xi\(\Phi(t)x\)$ are continuous in $t$, we have
$$\ba{ll}L(\Phi(s)x)&=\xi(\Phi(s)x)+\int_0^\8e^{-t}\xi(\Phi(t)\Phi(s)x)dt\\[1ex]
&\leq \xi(x)+\int_0^\8e^{-t}\xi(\Phi(t)\Phi(s)x)dt\\[1ex]
&<\xi(x)+\int_0^\8e^{-t}\xi(\Phi(t)x)dt=L(x).\ea$$

The other properties of $L$ simply  follow from that of $\xi$.  We omit the details of the argument.

If   $X$ is a normal Hausdorff space, then for any closed subset $K$ of $X$ with   $K\subset \W\- \sA$, there exists a nonnegative continuous  function $\psi$ on $X$ such that
$$
\psi(x)\equiv 0\,\,(\mb{on }\sA),\hs\mb{and } \psi(x)\equiv 1\,\,(\mb{on }K).
$$
Set $\eta(x)=\zeta(x)+\psi(x)$. Clearly $\zeta(x)$ is  a $\cK_0$ function of $\sA$  on $\W$.
If we replace the function $\zeta$ by $\eta$ in the above argument, one immediately obtains a Lyapunov function of $\sA$ satisfying
\eqref{cl}.   \,$\Box$
\br To guarantee the existence of a Lyapunov function of an attractor, we have assumed  in Theorem \ref{p2.13}  the existence of a $\cK_0$ function. In the case of a general topological space,  we do not know
whether such a function does exit. However, in many cases the space on which we are working may be the quotient space of a  pair $(N,E)$ of closed subsets of a metric space $X$. For such a space a $\cK_0$
 function can be directly formulated by using the metric $d$ of  $X$; see, e.g., \cite{LiLK}.\er

\section{ Morse Decompositions}
In this section we always assume $X$ is  a normal Hausdorff space.

 Let $\Phi$
be a given local semiflow on $X$, and $\sA$  an attractor of $\Phi$. Since $\sA$ is invariant, the restriction $\Phi|_\sA$ of $\Phi$ is a global semiflow on $\sA$.

A set  $A\subset \sA$ is called an {\em attractor} of $\Phi$ in $\sA$,
this means that  it is an attractor of $\Phi|_\sA$. We have the following fundamental result which generalizes the corresponding one for semiflows on metric spaces \cite{Kap}.
\bt Suppose $\sA$ has an admissible neighborhood $N_1$. Then any attractor $A$ in $\sA$ is also an attractor in $X$.
\et
\noindent{\bf Proof.} By definition $\sA$ attracts a neighborhood $N_2$ of $\sA$. As $X$ is normal, one can find a closed neighborhood $N$ of $\sA$ with $$N\subset N_1\cap N_2.$$ $N$ is strongly admissible. By stability of  $\sA$ there exists a closed neighborhood $M$ of $\sA$ such that
\be\label{e6.0}
\Phi(\R^+)M\subset N.
\ee

 Pick a closed neighborhood $\cO'$ of $A$ in $\sA$ with $\cO'\subset M$ such that $A$ attracts $\cO'$ (hence $\cO'\subset \W_\sA(A)$, where $\W_\sA(A)$ is the region of attraction of $A$ in $\sA$).
By the basic knowledge on general topology, one can find a closed neighborhood $\cO$ of $A$ in $X$ such that $$\cO'=\cO\cap \sA.$$ We can also assume $\cO\subset M$. We show that $A$ attracts a neighborhood $W\subset \cO$ of $A$, from which one immediately concludes that  $A$ is an attractor in $X$.

Let us first check that there is a neighborhood $W\subset \cO$ of $A$ such that \be\label{e6.1}\Phi(\R^+)W\subset\cO.\ee
Let $V=\mb{int}\cO$.  Then $V$ is an open neighborhood of $\ol A$.  For each integer $n>0$, since $\Phi(\R^+)x\subset \ol A\subset V$ for any fixed $x\in \ol A$ (see  Pro.\,\ref{p3.4}), by Lemma \ref{l2.3} there exists a neighborhood $U_x=U_x(n)$ of $x$ such that
$
\Phi([0,n])U_x\subset V.
$
Set $U_n=\Cup_{x\in \ol A}U_x$. Clearly $U_n$ is a neighborhood of $A$ for each $n$, and
\be\label{e6.2}
\Phi([0,n])U_n\subset V=\mb{int}\cO.\ee
We show that \be\label{e6.12}\Phi(\R^+)U_n\subset \cO\ee for some $n>0$, thus proving \eqref{e6.1}.

Suppose the contrary. Then for each $n$ there exists $x_n\in U_n$ such that $\Phi(t)x_n\not\in \cO$ for some $t>0$. Let
$$
t_n=\min\{t>0:\,\,\Phi(t)x_n\not\in \cO\}.
$$
Then \be\label{e6.3}\Phi(t_n)x_n\in \pa\cO,\hs \Phi([0,t_n])x_n\subset \cO.\ee
\eqref{e6.2} and the first relation in \eqref{e6.3} imply that $t_n>n$. Thus by admissibility of $\cO$ we deduce that the sequence $\Phi(t_n)x_n$ has a subsequence (still denoted by $\Phi(t_n)x_n$) converging to a point $x\in\pa\cO$.
For each $n$, let
$$
\gamma_n(t)=\Phi(t_n+t)x_n, \Hs t\in[-t_n,\8).
$$
Recalling that $\cO\subset M$,  by \eqref{e6.0} we have $$\gam_n(t)\in N,\Hs t\in [-t_n,\8)$$
for all $n$. Thanks to Pro.\,\ref{l2.9}, the sequence $\gam_n$ has a subsequence converging uniformly on any compact interval (in the sense given in Pro.\,\eqref{l2.9}) to a complete trajectory $\gam$. It is obvious that $\gam(0)=x\in\pa\cO$. Furthermore,
\be\label{e6.5}
\gam(t)\in \cO,\Hs t\in(-\8,0].
\ee
Because $\cO\subset M$, we infer from \eqref{e6.0} and \eqref{e6.5} that $\gam$ is contained in $N$. Thus by Theorem\,\ref{p2.3} (3) one concludes that $\gam(\R)\subset\sA$. \eqref{e6.5} then asserts that
\be\label{e6.6}
\gam(t)\in \cO\cap \sA=\cO'\subset \W_\sA(A),\Hs t\in(-\8,0].
\ee
Hence we actually have $\gam(\R)\subset\W_\sA(A)$.

In the following we further prove that $\gam(\R)\subset A$. For this purpose, we first verify that $\alpha(\gam)\subset A$.

By \eqref{e6.6} we see that $\alpha(\gam)\subset \cO'$. If $\alpha(\gam)\-A\ne\emp$,
then there is  a $y\in \alpha(\gam)\-A$. Pick a neighborhood $F$ of $A$ such that $y\not\in F$. As $A$ attracts $\cO'$, there exists $T>0$ such that
\be\label{e6.6b}
\Phi(t)\cO'\subset F,\Hs\A t>T.
\ee
On the other hand, by invariance of $\alpha(\gam)$, for any $t>T$ there exists $z\in \alpha(\gam)$ such that $\Phi(t)z=y$. Hence by  \eqref{e6.6b} one deduces that  $y\in F$, which leads to a contradiction! Thus $\alpha(\gam)\subset A$.

 Now we show that $\gam(\R)\subset \cO'$. It then follows by Theorem \ref{p2.3} (3) that $\gam(\R)\subset A$, which contradicts the fact that $\gam(0)=x\in\pa\cO$ and  completes the proof of what we desired in \eqref{e6.12}.

 By stability of $A$, one can find  a neighborhood $H$ of $A$ in $\sA$ with $H\subset \cO'$ such that \be\label{e6.11}\Phi(\R^+)H\subset \cO'.\ee Since $\alpha(\gam)\subset A$, we deduce that there is a $\tau<0$ such that $\gam(t)\in H$ for all $t<\tau$. Fix a $t_0<\tau$. Then by \eqref{e6.11} one has $\gam(t)\in \cO'$ for all $t\geq t_0$. Hence $\gam(\R)\subset \cO'$.

\vs
Now let $W$ be the neighborhood of $A$ in \eqref{e6.1}. Then $\w(W)\subset\cO$. Therefore we deduce by  invariance of $\w(W)$ and Theorem \ref{p2.3} (3) that $\w(W)\subset \sA$. Hence  $$\w(W)\subset \cO\cap \sA=\cO'.$$ Recalling that $\cO'\subset\W_\sA(A)$, again by Theorem \ref{p2.3} (3) one concludes that $\w(W)\subset A$. Thereby  $A$ attracts $W$.  $\Box$

\Vs

Let $A$ be an attractor of $\Phi$ in $\sA$. Define \be\label{e:3.0}
A^*=\sA\-\W_\sA(A). \ee
It is trivial to check that  $A^*$ is a nonempty $s$-compact invariant set. $A^*$ is called the {\em repeller }of
$\Phi$ in $\sA$  dual to $A$, and $(A,A^*)$  an {\em
attractor-repeller pair} in $\sA$.





\begin{proposition}\label{p6.4} Let $A$ be an attractor of $\Phi$ in
$\sA$, and let $\gamma:\R\ra \sA$ be a complete trajectory through $x\in
\sA$. Then the following properties hold. \vskip4pt $(1)$ If
$\,\w(\gamma)\cap A^*\neq\emp$, then $\gamma(\R)\subset A^*$.\vs
$(2)$ If $\,\alpha(\gamma)\cap \ol A^\sA\neq\emp$,  then $\gamma(\R)\subset
A$. Here $\ol A^\sA$ is the closure of $A$ in $\sA$.
\vs
$(3)$ If $x\in\sA\-( A\cup A^*)$, then $$\alpha(\gamma)\subset A^*,\hs \w(\gamma)\subset A.$$

\end{proposition}

\noindent{\bf Proof.} (1) Suppose $\gam(\tau)\not\in A^*$ for some $\tau\in\R$. Then by the definition of $A^*$ we have $\gam(\tau)\in \W_\sA(A)$. Since $A$ attracts $\gam(\tau)$, we have
 $\w(\gam)\subset \ol A^\sA$. Take a closed neighborhood $W$ of $A$ in $\sA$ with $W\subset \W_\sA(A)$.  Then by $s$-compactness of $\sA$ we see that $W$ is admissible. Thus  by Theorem \ref{p2.3} (3) one   deduces that
$\w(\gam)\subset A$. But this leads to  a contradiction. Hence $\gamma(\R)\subset A^*$.
\vs
(2) We first show that $\gamma(\R)\subset
\ol A^\sA$. Suppose the contrary. Then  $y=\gam(\tau)\not\in \ol A^\sA$ for some $\tau\in\R$. Take a neighborhood $V$ of $\ol A^\sA$ such that $y\not\in V$. By stability of $A$ in $\sA$ there is a neighborhood $W$ of $A$ in $\sA$ such that
$\Phi(\R^+)W\subset V$. On the other hand, since $\,\alpha(\gamma)\cap \ol A^\sA\neq\emp$, one can easily  find a $t_0<\tau$ such that $\gam(t_0)\in W$. It then follows that $\gam(t)\in V$ for all $t\geq t_0$. In particular,
$y=\gam(\tau)\in V$, a  contradiction! 

Now observe that $M=\gam(\R)\cup A$ is an invariant set in $\ol A^\sA$. Using a similar argument as in (1) one immediately concludes that $M\subset A$. Thereby $\gam(\R)\subset A$.
\vs
(3) As $x\not\in A^*$, we have $x\in \W_\sA(A)$. Thus the same argument in  (1) applies to show  that $\w(\gam)\subset A$.

In the following we verify that $\alpha(\gam)\subset A^*$. First, we infer from (2) that $\alpha(\gam)\cap\ol A^\sA= \emp$. Thus if $\alpha(\gam)\not\subset A^*$, then there is a point $y\in \alpha(\gam)$ such that
$y\in\W_\sA(A)\-\ol A^\sA$. Pick a neighborhood $U$ of $y$ in $\sA$ with $U\subset \W_\sA(A)$ and  a neighborhood $V$ of $A$ in $\sA$ such that  $U\cap V=\emp$. Then by attraction property of $A$   there exist a neighborhood $O$ of $y$ in $\sA$ with $O\subset U$ and a $T>0$ such that
\be\label{e6.9}
\Phi(t)O\subset V,\Hs t>T.
\ee
On the other hand, we infer from the definition of $\alpha$-limit set that there is a sequence $0<t_n\ra\8$ such that $\gam(-t_n)\in O$ for all $n$. Noticing that $y=\Phi(t_n)\gam(-t_n)$, by \eqref{e6.9} one finds that
$y\in V$, which leads to  a contradiction. \,$\Box$

\begin{definition}Let $\sA$ be an attractor.
An ordered collection $\cM=\{M_1,\cdots,M_n\}$ of subsets of $\sA$
is called a Morse decomposition of $\sA$, if there is an
increasing sequence
$
\emp=A_0\subset A_1\subset\cdots\subset A_n=\sA
$
of attractors of $\Phi$ in $\sA$ such that
\be\label{e6.10}
M_k=A_k\cap A_{k-1}^*,\Hs 1\leq k\leq n.
\ee

\end{definition}

The sets ${M_k}$ ($1\leq k\leq n$)  in  \eqref{e6.10}  will be referred to as  {\em Morse
sets} of $\sA$.

\begin{theorem}\label{t:3.2}

Let $\cM=\{M_1,\cdots,M_n\}$ be a Morse decomposition of $\sA$ with
the corresponding attractor sequence  $ \emp=A_0\subset
A_1\subset\cdots\subset A_n=\sA$. Then the following assertions hold.
\begin{enumerate}
\item[$(1)$] For each $k$, $(A_{k-1},M_k)$ is an
 attractor-repeller pair in $A_k$.

\item[$(2)$] ${M_k}$ ($1\leq k\leq n$) are pair-wise disjoint $s$-compact invariant sets.
 \item[$(3)$] If
$\gamma$ is a complete trajectory, then either $\gamma(\R)\subset
M_k$ for some Morse set $M_k$, or else there are indices $i<j$ such
that $ \alpha(\gamma)\subset M_j$ and $\w(\gamma)\subset M_i.$
\item[$(4)$] The attractors ${A_k}'s$ are uniquely determined by the
Morse sets. Specifically,
$$
A_k=\Cup_{1\leq i\leq k}W^u(M_i),\Hs 1\leq k\leq n,
$$
where $ W^u(M_i)$ is the unstable manifold of $M_i$ which is defined as
$$\ba{ll}
\mb{$W^u(M_i)=\{x|$ there is a  trajectory $\gamma:\R\ra \sA$}\\[1ex]
                  \hspace{2.8cm} \mb{through $x$ with $\alpha(\gamma)\subset M_i$\}.}\ea$$

\end{enumerate}
\end{theorem}

\noindent{\bf Proof.}  The proof of the theorem is the same as that of
\cite{Ryba}, Chapter III, Theorem 1.7, and is thus omitted. \,$\Box$

\end{document}